\UseRawInputEncoding
\documentclass[final,11pt,a4paper,reqno]{amsart}
\usepackage{stmaryrd}
\usepackage[left]{showlabels}
\usepackage[english]{babel}
\usepackage{pst-grad} 
\usepackage{pst-plot} 
\usepackage{pstricks}
\usepackage{amsmath,amssymb,mathrsfs,amsthm, mathtools}
\usepackage{tikz} 
\usepackage{mathcomp,wasysym}  
\usepackage{graphicx,subfigure}  
\usepackage[all]{xy} \xyoption{arc} \xyoption{color}
\usepackage{epsfig}
\usepackage{cite}
\usepackage[a4paper,
left=2.5cm, right=2.5cm,
top=3cm, bottom=3cm]{geometry}
\usepackage[bookmarks]{hyperref}
\hypersetup{
    colorlinks=true,
    linkcolor=cyan,
    filecolor=magenta,      
    urlcolor=cyan,
    citecolor=blue,
    pdftitle={},
    pdfpagemode=FullScreen,
    }
\usepackage{pgfplots}

\newcommand{\norm}[1]{\left\| #1 \right\|}

\normalsize
\normalsize
\newtheorem{theorem}{Theorem}

\theoremstyle{definition}
\newtheorem{rem}{Remark}

\pgfplotsset{
    every axis/.append style={
            axis x line=middle,
            axis y line=middle,
            xlabel={$x$},
            ylabel={$y$},
            axis line style={->},
        },
    marya/.style={color=black,thick,mark=none},
    soldot/.style={color=black,only marks,mark=*},
    holdot/.style={color=black,fill=white,only marks,mark=*},
    grid style={dotted,gray},
}

\tikzset{>=stealth}

\title[]{Well-posedness and decay for a nonlinear propagation wave model in atmospheric flows}

\author[D. Alonso-Or\'{a}n]{Diego Alonso-Or\'{a}n}
\email{dalonsoo@ull.edu.es}
\address{Departamento de An\'{a}lisis Matem\'{a}tico y Instituto de Matem\'{a}ticas y Aplicaciones (IMAULL), Universidad de La Laguna C/. Astrof\'{i}sico Francisco S\'{a}nchez s/n, 38200 - La Laguna, Spain.}

\author[R. Granero-Belinch\'{o}n]{Rafael Granero-Belinch\'{o}n}
\email{rafael.granero@unican.es}
\address{Departamento  de  Matem\'aticas,  Estad\'istica  y  Computaci\'on,  Universidad  de Cantabria.  Avda.  Los  Castros  s/n,  Santander,  Spain.}

\begin{document}
\begin{abstract}
In this note, we provide two results concerning the global well-posedness and decay of solutions to an asymptotic model describing the nonlinear wave propagation in the troposphere, namely, the morning glory phenomenon. The proof of the first result combines a pointwise estimate together with some interpolation inequalities to close the energy estimates in Sobolev spaces. The second proof relies on suitable Wiener-like functional spaces.
\end{abstract}



\maketitle

\allowdisplaybreaks

\section{Introduction and main result}
Many fascinating events that constantly test our comprehension of the dynamic and intricate atmospheric processes on Earth characterize the field of atmospheric science. Of all these mysterious events, \textit{the morning glory} is the one weather phenomenon that has fascinated scientists for decades. Long, horizontal cloud bands that often spread across the sky to make a recognizable and arresting pattern are what define the morning glory phenomenon. Recently in \cite{Constantin-Johnson-22}, the authors derived from the general Navier-Stokes equation in rotating spherical coordinates a more tractable asymptotic nonlinear system that describe this wave propagation given by
\begin{subequations}\label{eq:1}
\begin{align}
u_t+uu_{x}+vu_{y}&=\mu \Delta u+\alpha u+\beta v+F, \ \mbox{ in } \Omega, \ t>0, \\
u_{x}+v_{y}&=0,  \ \mbox{ in } \Omega, \ t>0.
\end{align}
\end{subequations}
The components of the vector velocity field are denoted by $u,v$, $\alpha, \beta\in\mathbb{R}$ and $\mu\in (0,\infty)$ is the viscosity parameter. The force $F$ represents a thermodynamic forcing term, which comprises the heat sources driving the motion. Equation (\ref{eq:1}b) comprises the incompressibility condition of the flow. In (\ref{eq:1}) the spatial domain $\Omega$ is the two-dimensional channel domain
\[ \Omega=\{ (x,y) \text{ s.t. }  x\in \mathbb{T}, \  0 < y < 1\}, \]
and the time satisfies $t\in [0,T]$ for certain $0<T\leq\infty$. Moreover, equations (\ref{eq:1}) are subject to the following boundary conditions
\begin{subequations}\label{bc:eq1}
\begin{align}
u&=0, \ \partial \Omega, \ t>0, \\
v&=0, \ \mbox{ on } \{y=0\}, \ t>0.
\end{align}
\end{subequations}
The corresponding initial-value problem consists of the system \eqref{eq:1}, \eqref{bc:eq1} along with the initial conditions
\begin{equation}\label{initialdata}
u(x,0)=u_{0}(x)
\end{equation}
which is assumed to be smooth enough for the purposes of the work (cf. the statement of Theorem \ref{thm1} for a the precise assumptions on the initial data). 

As already mentioned, system \eqref{eq:1} was derived by Constantin \& Johnson \cite{Constantin-Johnson-22}.  Moreover, Constantin \& Johnson provided a number of exact solutions: breeze-like flows, bore-like flows as well as oscillatory-like solutions, cf. \cite[\S 6]{Constantin-Johnson-22}. The same named authors investigated the existence of travelling-wave solutions, cf. \cite{Constantin-Johnson-23}. Their analysis relies on studying a nonlinear second ordinary differential equation by means of a global analysis phase-space using Lyapunov functions.  Recently, Matioc \& Roberti in \cite{Matioc-Roberti-23}, using an abstract quasilinear parabolic evolution framework, showed the global existence of weak solutions to \eqref{eq:1} as well as the local existence of classical solutions. 

In order to present the main result of this work, it is convenient to rewrite \eqref{eq:1} by eliminating the variable $v$. Following the approach in \cite{Matioc-Roberti-23}, we find that integrating $v$ in (\ref{eq:1}b) from $0$ to $y$, we have that
\[ v(x,y)=\int_{0}^{y} v_{y}(x,\xi) \ d\xi = -\int_{0}^{y} u_{x}(x,\xi) \ d\xi:= -Tu(x,y), \]
where $Tf:\Omega\to\mathbb{R}$ is given by
\[  Tf(x,y):=	\int_{0}^{y} f_{x}(x,\xi) \ d\xi. \]
This idea of using the fundamental theorem of calculus to express $v$ in terms of $u$ is reminiscing of the viscous primitive equations of large scale ocean and atmosphere dynamics, cf. \cite{Cao-Titi-07}. Therefore, using this observation  system \eqref{eq:1} can be rewritten as 
\begin{subequations}\label{eq:2}
\begin{align}
u_t+uu_{x}-Tu u_{y}&=\mu \Delta u+\alpha u- \beta Tu+F, \ \mbox{ in } \Omega, \ t>0, \\
Tu&=\int_{0}^{y} u_{x}(x,\xi) \ d\xi.
\end{align}
\end{subequations}
For convenience, we will work with a $2\pi$-periodic version of the infinite channel, namely, 
\[ \Omega=\{ (x,y) :  x\in \mathbb{T}, \  0 < y < 1\}, \]
and take $\alpha\leq 0$, and no external forcing, i.e, $F\equiv 0$. Hence, \eqref{eq:2} becomes
\begin{subequations}
\begin{align}
\label{eq:new31} u_t+uu_{x}-Tu u_{y}&=\mu \Delta u+\alpha u- \beta Tu, \ \mbox{ in } \Omega, \ t>0, \\
\label{eq:new32} Tu &=\int_{0}^{y} u_{x}(x,\xi) \ d\xi, 
\end{align}
\end{subequations}
supplemented with 
\begin{equation}\label{bc:3}
 u=0, \ \partial \Omega, \ t>0, \ u(z,0)=u_{0}(z), \ z\in \Omega.
 \end{equation}

The main result of this work is to provide the decay of the $L^\infty$ norm for arbitrary initial data together with a global existence of classical solutions to \label{eq:3} under a smallness assumption on the initial data. More precisely, the result reads as follows
\begin{theorem}\label{thm1}
Let $u_0\in H^{1}(\Omega)$ a zero mean function. Then, the corresponding solution satisfies
$$
\|u(t)\|_{L^\infty(\Omega)}\leq \|u_0\|_{L^\infty(\Omega)}. 
$$
Furthermore, if $u_0$ is such that $\norm{u_{0}}_{L^{\infty}(\Omega)}$ is sufficiently small. Then, there exists a unique global in time solution $u\in C([0,\infty),H^1(\Omega))\cap L^2(0,\infty;H^2(\Omega))$.
\end{theorem}
\begin{rem}
Just after the completion of this article, the preprint \cite{Host} appeared in the arXivs. The author shows a similar result to the one provided in this article but with some key differences. Compared with \cite{Host}, our work proves the $L^{\infty}$ decay of the solution. Such decay is not available in the literature before and to the best of the author's knowledge it is new. Furthermore, our well-posedness result requires an $L^{\infty}$ smallness assumption on the initial data to absorb the non-linear contribution. In \cite{Host} such nonlinear term is handled using by means of a clever and elegant extension method. 
\end{rem}
We also prove the following theorem
\begin{theorem}\label{thm2}
Let $u_0\in \tilde{A^0}(\Omega)$ a zero mean function such that $\norm{u_{0}}_{\tilde{A^0}(\Omega)}$ is sufficiently small. Then, there exists a unique global in time solution $u\in C([0,\infty),\tilde{A^0}(\Omega))\cap L^1(0,\infty; \tilde{A^2}(\Omega))$.
\end{theorem}

\subsection*{Notation} 
We denote by $(f,g)$ the $L^{2}$ inner product of $f,g$. For $m\in\mathbb{N}$, the natural inhomogeneous and homogeneous Sobolev norms  are defined by
\[ \norm{f}^{2}_{H^{m}(\Omega)}:= \norm{f}^{2}_{L^{2}(\Omega)} + \norm{\partial^{m}f}^{2}_{L^{2}(\Omega)}, \quad \norm{f}_{\dot{H}^{m}(\Omega)}^{2}:= \norm{\partial^{m}f}^{2}_{L^{2}(\Omega)}, \]
respectively. For notational convenience, we will just write $L^{2}, H^{m}, \dot{H}^{m}$ instead of $L^{2}(\Omega), H^{m}(\Omega), \dot{H}^{m}(\Omega).$ Moreover,  throughout the paper $C = C(\cdot)$ will denote a positive constant that may depend on fixed parameters  and  $x \lesssim y$ ($x \gtrsim y$) means that $x\le C y$ ($x\ge C y$) holds for some $C$.

\section{Proof of Theorem \ref{thm1}}\label{sec:proof:thm}
The proof follows by deriving appropriate a priori energy estimates combined with a suitable approximation procedure. To that purpose, we first provide the a priori energy estimates.
First of all, notice that mean-zero property is conserved as long as the solution exists.  Denoting by $z=(x,y)\in\Omega$ and integrating 
(\ref{eq:new31}) in $\Omega$ we find that
\begin{equation}\label{mean:zero:eq}
\partial_{t} \int_{\Omega} u(z,t) dz  = -\int_{\Omega} \left(uu_{x}-Tu u_{y}-\mu \Delta u-\alpha u+\beta Tu\right)(z,t) \ dz.
 \end{equation}
Notice that the first and the third term on the right hand side in \eqref{mean:zero:eq} are zero using the periodicity in the $x$ variable and the boundary condition \eqref{bc:3}. Furthermore, integrating by parts in the $y$ variable, recalling the definition (\ref{eq:new32}) and using that $u=0$ in $\partial \Omega$ we also have that
\[ \int_{\Omega} Tu u_{y} \ dz = -\int_{\Omega} u_{x} u \ dz=0. \]
Similarly, 
\[ \beta \int_{\Omega} Tu \ dz = -\int_{\Omega}\int_{0}^{y} u_{x}(x,\xi) d\xi \ dz=0. \]
Therefore,
\[ \partial_{t} \int_{\Omega} u(z,t) dz  =  \alpha \int_{\Omega} u(z,t) \ dz, \]
and since by assumption $u_{0}(z)$ has zero mean and $\alpha\leq0$ we conclude that 
\begin{equation}\label{mean:zero:evol}
 \int_{\Omega} u(z,t) \ dz=  0, \quad \mbox{for } t>0.
 \end{equation}
Next, let us derive a pointwise estimate. Following, \cite{Cordoba-Cordoba-04, CGO}, we define
\begin{align}
M(t)=\displaystyle\max_{z\in\Omega}u(z,t)=u(\overline{z}_{t},t), \mbox{ for } t>0, \\
m(t)=\displaystyle\min_{z\in\Omega}u(z,t)=u(\underline{z}_{t},t), \mbox{ for } t>0.
\end{align}
One can readily check that $M(t), m(t)$ are Lipschitz functions. Moreover, one can readily check that $M(t)$ satisfies
\begin{align*}
\left|M(t)-M(s)\right|&=\left\{\begin{array}{cc}u(\overline{z}_t,t)-u(\overline{z}_s,s) \text{ if } M(t)>M(s)\\
u(\overline{z}_s,s)-u(\overline{z}_t,t) \text{ if }M(s)>M(t)\end{array}\right.
\\
&\leq \left\{\begin{array}{cc}u(\overline{z}_t,t)-u(\overline{z}_t,s) \text{ if } M(t)>M(s)\\
u(\overline{z}_s,s)-u(\overline{z}_s,t) \text{ if }M(s)>M(t)\end{array}\right.
\\
&\leq \left\{\begin{array}{cc}|\partial_t u(\overline{z}_t,\xi)||t-s| \text{ if } M(t)>M(s)\\
|\partial_t u(\overline{z}_s,\xi)||t-s| \text{ if }M(s)>M(t)\end{array}\right.
\\
&\leq \max_{\eta,\xi}|\partial_t u(\eta,\xi)||t-s|.
\end{align*}
Similarly
$$
\left|m(t)-m(s)\right|\leq \max_{\eta,\xi}|\partial_t u(\eta,\xi)||t-s|.
$$
From Rademacher's theorem it holds that $M(t)$ and $m(t)$ are differentiable in $t$ almost everywhere. Furthermore
\begin{align*}
M'(t)&=\lim_{\delta\rightarrow0} \frac{u(\overline{z}_{t+\delta},t+\delta)-u(\overline{z}_t,t)}{\delta}=\lim_{\delta\rightarrow0} \frac{u(\overline{z}_{t+\delta},t+\delta)-u(\overline{z}_t,t)\pm u(\overline{z}_{t+\delta},t)}{\delta}\\
&\leq\lim_{\delta\rightarrow0} \frac{u(\overline{z}_{t+\delta},t+\delta)-u(\overline{z}_{t+\delta},t)}{\delta}\leq \partial_t u(\overline{z}_{t},t).
\end{align*}
In a similar fashion, we obtain that
\begin{align*}
M'(t)&=\lim_{\delta\rightarrow0} \frac{u(\overline{z}_{t+\delta},t+\delta)-u(\overline{z}_t,t)}{\delta}=\lim_{\delta\rightarrow0} \frac{u(\overline{z}_{t+\delta},t+\delta)-u(\overline{z}_t,t)\pm u(\overline{z}_{t},t+\delta)}{\delta}\\
&\geq\lim_{\delta\rightarrow0} \frac{u(\overline{z}_{t},t+\delta)-u(\overline{z}_{t},t)}{\delta}\geq \partial_t u(\overline{z}_{t},t).
\end{align*}
As a consequence
\begin{equation}\label{adg19}
M'(t)=\partial_t u(\overline{z}_{t},t) \text{ a.e.}
\end{equation}
Similarly
$$
m'(t)=\partial_t u(\underline{z}_{t},t) \text{ a.e.}
$$
Therefore, using (\ref{eq:new31}) and noticing that $u_{x}(\overline{z}_{t},t)=u_{y}(\overline{z}_{t},t)=0$ and $\Delta u(\overline{z}_{t},t) \leq 0$ we find that
\[ M'(t)\leq 0, \]
which implies that
\begin{equation}\label{est:max1}
 M(t)\leq M(0). 
 \end{equation}
Similarly, repeating the same argument and recalling that $\Delta u(\underline{z}_{t},t) \geq 0$,
\begin{equation}\label{est:min1}
 m(t)\geq m(0).
 \end{equation}

Notice that the maximum is obtained in the interior of $\Omega$, i.e. $z\in \dot{\Omega}$ and moreover $M(t)>0$. Otherwise, the maximum is obtained on the boundary and hence using the boundary condition \eqref{bc:3} this implies $M(t)=0$. Similarly, the minimum must be obtained in the interior of $\Omega$, otherwise $m(t)=0$ but this violates \eqref{mean:zero:evol}. As a consequence $m(t)<0$. Hence, combining \eqref{est:max1} and \eqref{est:min1} we have that
\begin{equation}\label{decrease:Linf}
\norm{u}_{L^{\infty}(\Omega)}\leq \norm{u_{0}}_{L^{\infty}(\Omega)}.
\end{equation}

Testing equation  (\ref{eq:new31}) against $u$, integrating by parts and using the periodicity in $x$ variable and boundary condition \ \eqref{bc:3} we have that
\begin{equation}\label{eq:L2:1}
\frac{1}{2}\frac{d}{dt} \norm{u}^{2}_{L^{2}} +\int_{\Omega} |\nabla u(z,t) |^{2} \ dz=-\int_{\Omega} Tu u_{y} u \ dz+\alpha\norm{u}^{2}_{L^2}-\beta\int_{\Omega} Tu u \ dz
\end{equation}
 Moreover, integrating by parts in the $y$ variable, recalling the definition (\ref{eq:new32}) and using that $u=0$ in $\partial \Omega$ we also have that 
 \begin{equation}\label{eq:L2:2}
  \int_{\Omega} Tu u_{y} u \ dz=- \frac{1}{2} \int_{\mathbb{T}}\int_{0}^{1} u_{x}(x,y) u^{2}(x,y) dx \ dy=0, 
  \end{equation}
by periodicity in the $x$ variable. Furthermore, by using Jensen's inequality one can readily check that 
\[ \norm{Tu}_{L^2}\leq \left( \int_{\mathbb{T}}\int_{0}^{1} \left( \int_{0}^{y} u_{x}(x,\xi) \ d\xi \right)^{2} dz \right)^{1/2} \leq \norm{\nabla u}_{L^{2}}, \]
which combined with H\"older's and Young's inequality shows that 
\begin{equation}\label{estimada:Tu}
\beta \int_{\Omega} |Tu u| \ dz \leq \beta \norm{Tu}_{L^{2}}\norm{u}_{L^{2}}\leq \frac{1}{2} \norm{\nabla u}_{L^{2}}^{2}+\frac{\beta^{2}}{2} \norm{u}^{2}_{L^{2}}. 
\end{equation}
Hence, using \eqref{eq:L2:1}, \eqref{eq:L2:2}, \eqref{estimada:Tu} and recalling that $\alpha\leq 0$ we have that
 \[\frac{d}{dt} \norm{u}^{2}_{L^{2}} +\int_{\Omega} |\nabla u(z,t) |^{2} \ dz\leq \beta^{2}  \norm{u}^{2}_{L^{2}},\]
 and thus Gr\"onwalls ineqality yields
 \begin{equation}\label{estimateL2}
  \displaystyle\sup_{t\in [0,T]}\norm{u(t)}^{2}_{L^{2}}+ \int_{0}^{T} \norm{\nabla u(t)}^{2}_{L^2} \ dt \leq \norm{u_{0}}^{2}_{L^{2}} e^{T}.
  \end{equation}
Next we compute the evolution of the $\dot{H}^{1}$ norm. To that purpose, testing (\ref{eq:new31}) against $-\Delta u$ and integrating by parts
\begin{align*}
\frac{1}{2}\frac{d}{dt} \norm{\nabla u}^{2}_{L^2}+ \int_{\Omega} |\Delta u|^{2} \ dz&=\int_{\Omega} Tu u_{y} \Delta u \ dz-\int_{\Omega} u u_{x} \Delta u \ dz  -\alpha\int_{\Omega} u \Delta u \ dz+\beta\int_{\Omega} Tu \Delta u  dz   \\
&=I_{1}+I_{2}+I_{3}+I_{4}.
\end{align*}
By means of H\"older's and Jensen's inequality we find that
\begin{align}\label{i1i2}
|I_{1}+I_{2}| &\leq \norm{Tu}_{L^{4}}\norm{u_{y}}_{L^{4}}\norm{\Delta u}_{L^{2}}+\norm{u}_{L^{\infty}}\norm{u_{x}}_{L^{2}}\norm{\Delta u}_{L^2} \nonumber \\
&\leq  \norm{u_{x}}_{L^{4}}\norm{u_{y}}_{L^{4}}\norm{\Delta u}_{L^{2}}+\norm{u}_{L^{\infty}}\norm{u_{x}}_{L^{2}}\norm{\Delta u}_{L^2}. 
\end{align}
Furthermore, integrating by parts in $x$ and recalling the periodic boundary conditions in $x$, we have that
\begin{equation} \label{interpolation:1}
\norm{u_{x}}^{4}_{L^{4}}\leq 3 \norm{u_{x}}^{2}_{L^{4}}\norm{u_{xx}}_{L^{2}}\norm{u}_{L^{\infty}} \Longleftrightarrow \norm{u_{x}}^{2}_{L^{4}}\leq 3 \norm{u_{xx}}_{L^{2}}\norm{u}_{L^{\infty}}.
\end{equation}
Similarly, using the fact that $u=0$ on $\partial\Omega$, we infer that
\begin{equation}\label{interpolation:2}
\norm{u_{y}}^{2}_{L^{4}}\leq 3 \norm{u_{yy}}_{L^{2}}\norm{u}_{L^{\infty}}.
\end{equation}
Hence, using \eqref{interpolation:1}-\eqref{interpolation:2} into \eqref{i1i2} and applying Young's inequality we arrive to
\begin{align}
|I_{1}+I_{2}|&\leq  3 \norm{u}_{L^{\infty}}\norm{\Delta u}^{2}_{L^{2}}+\norm{u}_{L^{\infty}}\norm{u_{x}}_{L^{2}}\norm{\Delta u}_{L^2} \nonumber \\
&\leq   
\left( 3\norm{u}_{L^{\infty}}+\frac{1}{2}\norm{u}^{2}_{L^{\infty}}\right)\norm{\Delta u}^{2}_{L^{2}}+ \frac{1}{2}\norm{\nabla u}_{L^2}^{2}.\label{i1i2:final}
\end{align}
The remaining terms $I_{3}, I_{4}$ can be estimated easily, namely,
\begin{align}
|I_{3}|&\leq |\alpha| \norm{\Delta u}_{L^{2}}\norm{u}_{L^{2}} \leq \frac{\epsilon}{2} \norm{\Delta u}_{L^{2}}^{2} + \frac{|\alpha|^{2}}{2\epsilon}\norm{u}_{L^{2}}^{2}, \label{i3:final} \\
|I_{4}| &\leq  C  \norm{Tu}_{L^{2}}\norm{\Delta u}_{L^{2}}\leq \frac{\epsilon}{2} \norm{\Delta u}_{L^{2}}^{2} + \frac{|\beta|^{2}}{2\epsilon}\norm{\nabla u}_{L^{2}}^{2}. \label{i4:final}
\end{align}

Hence, collecting bounds \eqref{i1i2:final}, \eqref{i3:final}, \eqref{i4:final}, we infer that
\begin{equation}
\frac{1}{2}\frac{d}{dt} \norm{\nabla u}^{2}_{L^2}+ \int_{\Omega} |\Delta u|^{2} \ dz \leq \left( 3\norm{u}_{L^{\infty}}+\frac{1}{2}\norm{u}^{2}_{L^{\infty}}+\epsilon\right)\norm{\Delta u}^{2}_{L^{2}}+ \left(\frac{1}{2}+\frac{|\beta|^{2}}{2\epsilon}\right)\norm{\nabla u}_{L^2}^{2}+ \frac{|\alpha|^{2}}{2\epsilon}\norm{u}_{L^{2}}^{2}.
\end{equation}
Invoking the $L^{\infty}$ decay estimate \eqref{decrease:Linf} and taking $\epsilon$ and $\norm{u_{0}}_{L^{\infty}}$ small enough we conclude that
\begin{equation}
\frac{1}{2}\frac{d}{dt} \norm{\nabla u}^{2}_{L^2}+ \frac{1}{2}\int_{\Omega} |\Delta u|^{2} \ dz \leq C_{\beta} \norm{\nabla u}_{L^2}^{2}+ C_{\alpha} \norm{u}_{L^{2}}^{2}.
\end{equation}
Integrating in time and taking supremums, we have that
\begin{align}
\displaystyle\sup_{t\in[0,T]} \norm{\nabla u}^{2}_{L^2}+ \int_{0}^{T}\int_{\Omega} |\Delta u|^{2} \ dz &\leq C_{\beta} \int_{0}^{T} \norm{\nabla u}_{L^2}^{2} \ dt+ C_{\alpha} \int_{0}^{T} \norm{u}_{L^{2}}^{2} \ dt  \nonumber \\
&\leq \norm{\nabla u_{0}}^{2}_{L^{2}} + C_{\beta} \norm{u_{0}}_{L^{2}}e^{T} + C_{\alpha} \norm{u_{0}}_{L^{2}}e^{T} T
\end{align}
where in the second inequality we have used estimate \eqref{estimateL2}. This concludes the \emph{a priori} estimates. The approximation procedure can be done using Galerkin method. To obtain the uniqueness, we argue by contradiction. If $u$ and $w$ are two solutions emanating from the same initial data, we consider $U$ their difference. Then, $U$ solves
$$
U_t=-uU_{x}-TU u_{y}-Uw_{x}-Tw U_{y}+\mu \Delta U+\alpha U- \beta TU.
$$
Testing against $U$ and integrating by parts, we find
$$
\frac{d}{dt}\|U\|_{L^2}^2\leq -\|\nabla U\|_{L^2}^2+\|U\|_{L^2}\|\nabla U\|_{L^2}(\|u\|_{L^\infty}+\|w\|_{L^\infty})+\|U_y\|_{L^2}\|T U\|_{L^2}\|u\|_{L^\infty}+2\|U\|_{L^2}\|U_{x}\|_{L^2}\|w\|_{L^\infty}.
$$
Using the smallness hypothesis we conclude the desired bound and the result.

\section{Proof of Theorem \ref{thm2}}\label{sec:proof:thm2}
Since
\[ \Omega=\{ (x,y) \text{ s.t. }  x\in \mathbb{T}, \  0 < y < 1\}, \]
we consider the Fourier series representation
$$
u(x,y)=\sum_{n\in\mathbb{Z}}\sum_{m\geq1}\hat{u}(n,m)e^{inx}\sin(m\pi y),
$$
together with the Wiener-like spaces
$$
\tilde{\mathcal{A}^s}=\left\{u\text{ s.t. }\sum_{n\in\mathbb{Z}}\sum_{m\geq1}(|n|^s+|m|^s)|\hat{u}(n,m)|<\infty\right\}.
$$
These spaces have similar properties to those in \cite{BGB} and the references therein. In particular, we observe that 
$$
\tilde{\mathcal{A}^s}\subset C^s
$$
and that they form a Banach scale of Banach algebras. These spaces will allow us to achieve maximal parabolic regularity. Indeed, if we now estimate the evolution of these norms, we find that
$$
\frac{d}{dt}\|u(t)\|_{\tilde{\mathcal{A}^0}}\leq \|u\|_{\tilde{\mathcal{A}^0}}\|u\|_{\tilde{\mathcal{A}^1}}+\|u\|_{\tilde{\mathcal{A}^1}}^2-\mu \|u\|_{\tilde{\mathcal{A}^2}},
$$
where we have used $\alpha\leq 0$ and
$$
\widehat{Tu}(n,m) =\frac{n}{m}\hat{u}(n,m). 
$$
Using interpolation, we find that
$$
\|u\|_{\tilde{\mathcal{A}^1}}^2\leq \|u\|_{\tilde{\mathcal{A}^0}}\|u\|_{\tilde{\mathcal{A}^2}}.
$$
As a consequence, we conclude
$$
\frac{d}{dt}\|u(t)\|_{\tilde{\mathcal{A}^0}}\leq 2\|u\|_{\tilde{\mathcal{A}^0}}\|u\|_{\tilde{\mathcal{A}^2}}-\mu \|u\|_{\tilde{\mathcal{A}^2}}.
$$
This concludes the \emph{a priori} estimates of the solution. The approximation procedure can be done using a standard Galerkin approximation in finite dimensional projections. The uniqueness is a consequence of the obtained regularity and a standard contradiction result, from where we conclude the desired result.

\section*{Acknowledgments}
\noindent  D.A-O is supported by the Spanish MINECO through Juan de la Cierva fellowship FJC2020-046032-I.  Both authors are funded by the project "An\'alisis Matem\'atico Aplicado y Ecuaciones Diferenciales" Grant PID2022-141187NB-I00 funded by MCIN /AEI /10.13039/501100011033 / FEDER, UE and acronym "AMAED". This publication is part of the project PID2022-141187NB-I00 funded by MCIN/ AEI /10.13039/501100011033.

\begin{footnotesize}

\end{footnotesize}
\vspace{2cm}

\end{document}